\documentclass[11pt]{article}
\usepackage{amsmath}
\usepackage{amssymb}
\usepackage{amscd}
\usepackage{xspace}
\usepackage{verbatim}
\newcommand{\mysection}[1]{
\section{#1}\setcounter{equation}{0}}
\date{}
\begin{document}
\begin{center}{\bf \Large Existence and stability of solutions of general\\[2mm]

\noindent  semilinear elliptic equations with measure data}

 \end{center}
\begin{center}{\bf  Laurent V\'eron}\\
{\small Laboratoire de Math\'ematiques et Physique Th\'eorique}\\
 {\small  Universit\'e Fran\c cois-Rabelais, Tours,  FRANCE} \\[2mm]

 \end{center}



\newcommand{\txt}[1]{\;\text{ #1 }\;}
\newcommand{\tbf}{\textbf}
\newcommand{\tit}{\textit}
\newcommand{\tsc}{\textsc}
\newcommand{\trm}{\textrm}
\newcommand{\mbf}{\mathbf}
\newcommand{\mrm}{\mathrm}
\newcommand{\bsym}{\boldsymbol}
\newcommand{\scs}{\scriptstyle}
\newcommand{\sss}{\scriptscriptstyle}
\newcommand{\txts}{\textstyle}
\newcommand{\dsps}{\displaystyle}
\newcommand{\fnz}{\footnotesize}
\newcommand{\scz}{\scriptsize}
\newcommand{\be}{
\begin{equation}
}
\newcommand{\bel}[1]{
\begin{equation}
\label{#1}}
\newcommand{\ee}{
\end{equation}
}
\newcommand{\eqnl}[2]{
\begin{equation}
\label{#1}{#2}
\end{equation}
}
\newtheorem{subn}{\name}
\newcommand{\bsn}[1]{\def\name{#1}
\begin{subn}}
\newcommand{\esn}{
\end{subn}}
\newtheorem{sub}{\name}[section]
\newcommand{\dn}[1]{\def\name{#1}}   
\newcommand{\bs}{
\begin{sub}}
\newcommand{\es}{
\end{sub}}
\newcommand{\bsl}[1]{
\begin{sub}\label{#1}}
\newcommand{\bth}[1]{\def\name{Theorem}
\begin{sub}\label{t:#1}}
\newcommand{\blemma}[1]{\def\name{Lemma}
\begin{sub}\label{l:#1}}
\newcommand{\bcor}[1]{\def\name{Corollary}
\begin{sub}\label{c:#1}}
\newcommand{\bdef}[1]{\def\name{Definition}
\begin{sub}\label{d:#1}}
\newcommand{\bprop}[1]{\def\name{Proposition}
\begin{sub}\label{p:#1}}
\newcommand{\R}{\eqref}
\newcommand{\rth}[1]{Theorem~\ref{t:#1}}
\newcommand{\rlemma}[1]{Lemma~\ref{l:#1}}
\newcommand{\rcor}[1]{Corollary~\ref{c:#1}}
\newcommand{\rdef}[1]{Definition~\ref{d:#1}}
\newcommand{\rprop}[1]{Proposition~\ref{p:#1}} 
\newcommand{\BA}{
\begin{array}}
\newcommand{\EA}{
\end{array}}
\newcommand{\BAN}{\renewcommand{\arraystretch}{1.2}
\setlength{\arraycolsep}{2pt}
\begin{array}}
\newcommand{\BAV}[2]{\renewcommand{\arraystretch}{#1}
\setlength{\arraycolsep}{#2}
\begin{array}}
\newcommand{\BSA}{
\begin{subarray}}
\newcommand{\ESA}{\end{subarray}}
\newcommand{\BAL}{\begin{aligned}}
\newcommand{\EAL}{\end{aligned}}
\newcommand{\BALG}{\begin{alignat}}
\newcommand{\EALG}{\end{alignat}}
\newcommand{\BALGN}{\begin{alignat*}}
\newcommand{\EALGN}{\end{alignat*}}
\newcommand{\note}[1]{\textit{#1.}\hspace{2mm}}
\newcommand{\Proof}{\note{Proof}}
\newcommand{\qeda}{\hspace{10mm}\hfill $\square$}
\newcommand{\qed}{\\
${}$ \hfill $\square$}
\newcommand{\Remark}{\note{Remark}}
\newcommand{\modin}{$\,$\\
[-4mm] \indent}
\newcommand{\forevery}{\quad \forall}
\newcommand{\set}[1]{\{#1\}}
\newcommand{\setdef}[2]{\{\,#1:\,#2\,\}}
\newcommand{\setm}[2]{\{\,#1\mid #2\,\}}
\newcommand{\lra}{\longrightarrow}
\newcommand{\sgn}{\rm{sgn}}
\newcommand{\lla}{\longleftarrow}
\newcommand{\llra}{\longleftrightarrow}
\newcommand{\Lra}{\Longrightarrow}
\newcommand{\Lla}{\Longleftarrow}
\newcommand{\Llra}{\Longleftrightarrow}
\newcommand{\warrow}{\rightharpoonup}
\newcommand{
\paran}[1]{\left (#1 \right )}
\newcommand{\sqbr}[1]{\left [#1 \right ]}
\newcommand{\curlybr}[1]{\left \{#1 \right \}}
\newcommand{\abs}[1]{\left |#1\right |}
\newcommand{\norm}[1]{\left \|#1\right \|}
\newcommand{
\paranb}[1]{\big (#1 \big )}
\newcommand{\lsqbrb}[1]{\big [#1 \big ]}
\newcommand{\lcurlybrb}[1]{\big \{#1 \big \}}
\newcommand{\absb}[1]{\big |#1\big |}
\newcommand{\normb}[1]{\big \|#1\big \|}
\newcommand{
\paranB}[1]{\Big (#1 \Big )}
\newcommand{\absB}[1]{\Big |#1\Big |}
\newcommand{\normB}[1]{\Big \|#1\Big \|}

\newcommand{\thkl}{\rule[-.5mm]{.3mm}{3mm}}
\newcommand{\thknorm}[1]{\thkl #1 \thkl\,}
\newcommand{\trinorm}[1]{|\!|\!| #1 |\!|\!|\,}
\newcommand{\bang}[1]{\langle #1 \rangle}
\def\angb<#1>{\langle #1 \rangle}
\newcommand{\vstrut}[1]{\rule{0mm}{#1}}
\newcommand{\rec}[1]{\frac{1}{#1}}
\newcommand{\opname}[1]{\mbox{\rm #1}\,}
\newcommand{\supp}{\opname{supp}}
\newcommand{\dist}{\opname{dist}}
\newcommand{\myfrac}[2]{{\displaystyle \frac{#1}{#2} }}
\newcommand{\myint}[2]{{\displaystyle \int_{#1}^{#2}}}
\newcommand{\mysum}[2]{{\displaystyle \sum_{#1}^{#2}}}
\newcommand {\dint}{{\displaystyle \int\!\!\int}}
\newcommand{\q}{\quad}
\newcommand{\qq}{\qquad}
\newcommand{\hsp}[1]{\hspace{#1mm}}
\newcommand{\vsp}[1]{\vspace{#1mm}}
\newcommand{\ity}{\infty}
\newcommand{\prt}{
\partial}
\newcommand{\sms}{\setminus}
\newcommand{\ems}{\emptyset}
\newcommand{\ti}{\times}
\newcommand{\pr}{^\prime}
\newcommand{\ppr}{^{\prime\prime}}
\newcommand{\tl}{\tilde}
\newcommand{\sbs}{\subset}
\newcommand{\sbeq}{\subseteq}
\newcommand{\nind}{\noindent}
\newcommand{\ind}{\indent}
\newcommand{\ovl}{\overline}
\newcommand{\unl}{\underline}
\newcommand{\nin}{\not\in}
\newcommand{\pfrac}[2]{\genfrac{(}{)}{}{}{#1}{#2}}

\def\ga{\alpha}     \def\gb{\beta}       \def\gg{\gamma}
\def\gc{\chi}       \def\gd{\delta}      \def\ge{\epsilon}
\def\gth{\theta}                         \def\vge{\varepsilon}
\def\gf{\phi}       \def\vgf{\varphi}    \def\gh{\eta}
\def\gi{\iota}      \def\gk{\kappa}      \def\gl{\lambda}
\def\gm{\mu}        \def\gn{\nu}         \def\gp{\pi}
\def\vgp{\varpi}    \def\gr{\rho}        \def\vgr{\varrho}
\def\gs{\sigma}     \def\vgs{\varsigma}  \def\gt{\tau}
\def\gu{\upsilon}   \def\gv{\vartheta}   \def\gw{\omega}
\def\gx{\xi}        \def\gy{\psi}        \def\gz{\zeta}
\def\Gg{\Gamma}     \def\Gd{\Delta}      \def\Gf{\Phi}
\def\Gth{\Theta}
\def\Gl{\Lambda}    \def\Gs{\Sigma}      \def\Gp{\Pi}
\def\Gw{\Omega}     \def\Gx{\Xi}         \def\Gy{\Psi}

\def\CS{{\mathcal S}}   \def\CM{{\mathcal M}}   \def\CN{{\mathcal N}}
\def\CR{{\mathcal R}}   \def\CO{{\mathcal O}}   \def\CP{{\mathcal P}}
\def\CA{{\mathcal A}}   \def\CB{{\mathcal B}}   \def\CC{{\mathcal C}}
\def\CD{{\mathcal D}}   \def\CE{{\mathcal E}}   \def\CF{{\mathcal F}}
\def\CG{{\mathcal G}}   \def\CH{{\mathcal H}}   \def\CI{{\mathcal I}}
\def\CJ{{\mathcal J}}   \def\CK{{\mathcal K}}   \def\CL{{\mathcal L}}
\def\CT{{\mathcal T}}   \def\CU{{\mathcal U}}   \def\CV{{\mathcal V}}
\def\CZ{{\mathcal Z}}   \def\CX{{\mathcal X}}   \def\CY{{\mathcal Y}}
\def\CW{{\mathcal W}} \def\CQ{{\mathcal Q}} 
\def\BBA {\mathbb A}   \def\BBb {\mathbb B}    \def\BBC {\mathbb C}
\def\BBD {\mathbb D}   \def\BBE {\mathbb E}    \def\BBF {\mathbb F}
\def\BBG {\mathbb G}   \def\BBH {\mathbb H}    \def\BBI {\mathbb I}
\def\BBJ {\mathbb J}   \def\BBK {\mathbb K}    \def\BBL {\mathbb L}
\def\BBM {\mathbb M}   \def\BBN {\mathbb N}    \def\BBO {\mathbb O}
\def\BBP {\mathbb P}   \def\BBR {\mathbb R}    \def\BBS {\mathbb S}
\def\BBT {\mathbb T}   \def\BBU {\mathbb U}    \def\BBV {\mathbb V}
\def\BBW {\mathbb W}   \def\BBX {\mathbb X}    \def\BBY {\mathbb Y}
\def\BBZ {\mathbb Z}

\def\GTA {\mathfrak A}   \def\GTB {\mathfrak B}    \def\GTC {\mathfrak C}
\def\GTD {\mathfrak D}   \def\GTE {\mathfrak E}    \def\GTF {\mathfrak F}
\def\GTG {\mathfrak G}   \def\GTH {\mathfrak H}    \def\GTI {\mathfrak I}
\def\GTJ {\mathfrak J}   \def\GTK {\mathfrak K}    \def\GTL {\mathfrak L}
\def\GTM {\mathfrak M}   \def\GTN {\mathfrak N}    \def\GTO {\mathfrak O}
\def\GTP {\mathfrak P}   \def\GTR {\mathfrak R}    \def\GTS {\mathfrak S}
\def\GTT {\mathfrak T}   \def\GTU {\mathfrak U}    \def\GTV {\mathfrak V}
\def\GTW {\mathfrak W}   \def\GTX {\mathfrak X}    \def\GTY {\mathfrak Y}
\def\GTZ {\mathfrak Z}   \def\GTQ {\mathfrak Q}

\font\Sym= msam10 
\def\SYM#1{\hbox{\Sym #1}}
\newcommand{\bdw}{\prt\Gw\xspace}
\medskip

\noindent{\small {\bf Abstract} We study existence and stability for solutions of $ -L u+g(x,u)=\gw$ where $L$ is a second order elliptic operator, $g$ a Caratheodory function and $\gw$ a measure in $\overline\Gw$. We present a unified theory of the Dirichlet problem and the Poisson equation. We prove the stability of the problem with respect to weak convergence of the data. 

\noindent
{\it \footnotesize 2010 Mathematics Subject Classification}. {\scriptsize
35J61, 35J66, 28A20}.\\
{\it \footnotesize Key words:} {\scriptsize Elliptic operators, Borel measures, Marcinkiewicz spaces, $\Gd_2$ condition.
}
}

\mysection{Introduction}
 Let $\Gw$ be a smooth bounded domain of $\BBR^N$, $L$ a uniformly elliptic second order differential operator in divergence form with Lipschitz continuous coefficients and $g$ is a real valued Caratheodory function defined in $\Gw\ti\BBR$. If $\gw$ is a Radon measure on $\overline\Gw$, we study existence and stability of solutions of the generalized equation
 \bel{I1}
 -L u+g(x,u)=\gw
\ee
in $\overline\Gw$. Precise assumptions are made on the coefficients of $L$ so that uniqueness holds. A fundamental contribution is made by Benilan and Brezis \cite{Br1}, \cite{BeBr} who study the case where $L=\Gd$ and $g:\BBR\mapsto\BBR$ is nondecreasing and positive on $\BBR_+$: if  $\gm$ is a bounded measure in $\Gw$ and $g$ satisfies the {\it subcriticality assumption}
 \bel{I2}
\myint{1}{\infty}\left(g(s)-g(-s)\right)s^{-2\frac{N-1}{N-2}}ds<\infty,
\ee
then there exists a unique function $u\in L^1(\Gw)$ such that $g\circ u\in  L^1(\Gw)$ (where $g\circ u(x)= g(x,u(x))$) satisfying 
\bel{I3}
\myint{\Gw}{}\left(-u\Gd \gz+g\circ u\,\gz\right)dx=\myint{\Gw}{}\gz d\gm,
\ee
for any $\gz\in C^2_0(\Gw)$.\smallskip

The boundary value problem with measures is first investigated by  Gmira and V\'eron \cite{GV}. By adapting the method introduced by Benilan and Brezis they obtain the existence and uniqueness of a weak solution of 
 \bel{I4}\BA {ll}
 -\Gd u+g(u)=0\qquad&\text{in }\Gw\\
 \phantom{ -\Gd +g(u)}
 u=\gl&\text{in }\prt\Gw
\EA\ee
when $\gl$ is a Radon measure. They assume that $g$, always nondecreasing, satisfies the {\it boundary subcriticality assumption}
 \bel{I5}
\myint{1}{\infty}\left(g(s)-g(-s)\right)s^{-\frac{2N}{N-2}}ds<\infty,
\ee
and prove the existence and uniqueness of a weak solution to $(\ref{I4})$. For this problem, in the integral identity $(\ref{I3})$ the right hand-side is replaced by $-\myint{\prt\Gw}{}\gz_{\bf n} d\gl$ (where $\gz_{\bf n}=\nabla u.{\bf n}$ is the outward normal derivative on $\prt\Gw$). \smallskip

In \cite{Ve1} V\'eron extends Benilan-Brezis  results in replacing $\Gd$ by a general uniformly elliptic second order differential operator with smooth coefficients. If $g$ is nondecreasing and satisfies, for some $\ga\in [0,1]$, the {\it $\ga$-subcriticality assumption},
 \bel{I7}
\myint{1}{\infty}\left(g(s)-g(-s)\right)s^{-2\frac{N+\ga-1}{N+\ga-2}}ds<\infty,
\ee
then if $\gm$ belongs to $\mathfrak M_{\gr^\ga}(\Gw)$, which means
 \bel{I6}
\norm{\gm}_{\mathfrak M_{\gr^\ga}}:=\myint{\overline\Gw}{}\gr^{\ga}d\abs\gm<\infty,
\ee
where $\gr(x):=\dist (x,\prt\Gw)$, there exists a unique $u\in L^1(\Gw)$ such that $g(u)\in L^1_\gr(\Gw)$ satisfying
\bel{I8}
\myint{\Gw}{}\left(-uL^* \gz+g(u)\gz\right)dx=\myint{\Gw}{}\gz d\gm\qquad\forall\gz\in C_{c}^{1,L^*}(\overline \Omega).
\ee
 where
\bel{I9}
C_{c}^{1,L^*}(\overline \Omega)=\{\gz\in C^1(\overline \Gw):\gz=0\text{ on }\prt\Gw,\;L^*\gz\in 
L^{\ity}(\Gw)\},
\ee
where $L^*$ is the adjoint operator to $L$.  Furthermore he proves the weak stability of the problem. it means that if $u_n$ is a set of solutions of 
 \bel{I1-n}\BA {ll}
 -L u_n+g(u_n)=\gm_n\qquad&\text{in }\Gw\\\phantom{ -L +g(u_n)}
 u_n=0\qquad&\text{in }\prt\Gw
\EA\ee
for a sequence of measure $\{\gm_n\}$ such that
 \bel{I10}
\lim_{n\to\infty}\myint{\Gw}{}\gz d\gm_n=\myint{\Gw}{}\gz d\gm
\ee
for all $\gz\in C(\overline\Gw)$ verifying $\sup_{\Gw}\gr^{-\ga}|\gz|<\infty$, then $u_n\to u$ where $u$ satisfies $(\ref{I1})$. 
However, a careful observation of the existence and stability statements proved in \cite[Th 3.7, Cor 3.8]{Ve1} shows that the result is slightly stronger than the one stated since it implies the following: \smallskip

\noindent{  \it  Let $\ga\in [0,1]$ and $g:\BBR\mapsto\BBR$ be continuous function which satisfies the $\ga$-subcriticality assumption $(\ref{I7})$. If $\{\gm_n\}$ is a sequence of Radon measures in $\overline\Gw$ such that 
 \bel{I11}
\myint{\overline\Gw}{}\gr^{\ga}d\abs{\gm_n}\leq M
\ee
for some $M>0$ and $(\ref{I10})$ holds for $\gz$ such that $\gr^{-\ga}\gz\in C(\overline\Gw)$, then the corresponding solution $u_n$ of $(\ref{I1-n})$ converges to the solution $u$ of $(\ref{I1})$. In particular, if $\ga=1$, it contains the case where there exists a Radon measure $\gl $ on $\prt\Gw$ such that}
 \bel{I11'}
\lim_{n\to\infty}\myint{\Gw}{}\gz d\gm_n=-\myint{\prt\Gw}{}\gz_{\bf n} d\gl\qquad\forall \gz\in C_c^1(\overline\Gw). 
\ee
\smallskip

The case where the nonlinearity $g$ depends on the $\gr(x)$ variable has investigated by Marcus \cite{Ma}. If $g(x,r)\rm {sign\,}r\leq \gr(x)^\gb \tilde g(\abs r)\rm {sign\,}r$ for some $\gb>-2$ and $\tilde g$ satisfying a subcriticality assumption
 \bel{I12}
\myint{1}{\infty}\left(\tilde g(s)-\tilde g(-s)\right)s^{-\frac{2N+\gb-1}{N-1}}ds<\infty,
\ee
then there exists a weak solution to problem $(\ref{I4})$ for any Radon measure $\gl$. Furthermore stability holds. \smallskip

The subcriticality is a key hypothesis in all the previous results: essentially it means that the problem can be solved for any measure if it can be solved for a Dirac measure. The different integral assumptions are just the transcription that the fact that $g$ of the fundamental solution of the associated linear equation is integrable for a suitable measure associated to the distance function $\gr$.  \medskip

The aim of this article is twofold: 1- to unify the problems for measures in $\Gw$ and on $\prt\Gw$; 2- to present under the form of an integrability condition a classical sufficient condition of solvability which has the advantage of being a natural extension to the 
supercritical case the previous subcriticality assumptions and to provide new results results of existence and stability for $(\ref{I1})$ in the spirit of \cite{Ve1}. A function $g:\Gw\ti\BBR\mapsto\BBR$ belongs to the class $G_{h,\Psi}$ if it is a Caratheodory function and there exist  a continuous and nondecreasing function $\tilde g:\BBR\mapsto\BBR$ vanishing at $0$,  a locally integrable nonnegative function $h$ defined in $\Gw$ and a nonnegative continuous nonincreasing function  $\Psi:[0,\infty)\mapsto[0,\infty)$,  such that
 \bel{I13}\BA {ll}
\abs{g(x,r)}\leq h(x) \abs{\tilde g(r)}\qquad\forall (x,r)\in \Gw\ti\BBR,
\EA\ee
and the {\it $\Psi$-integrability condition} holds, i.e. 
 \bel{I14}
-\myint{0}{\infty}\left(\tilde g(s)-\tilde g(-s)\right)d\Psi(t)ds<\infty.
\ee

Let $G$ and $K$ be respectively the Green and Poisson kernels corresponding to the operator $L$ in $\Gw$ and $\BBG[.]$ and $\BBK[.]$ the corresponding potential operators. The natural subcritical assumptions in the framework of Marcus's results (with $h$ instead of $\gr^\gb$) for solving 
 \bel{I15}\BA {ll}
-Lu+g(x,u)=\gm\qquad&\text{in }\Gw\\[2mm]\phantom{-L+G(x,u)}
u=\gl&\text{in }\prt\Gw
\EA\ee
would be
 \bel{I16}
\myint{1}{\infty}\left(\BBG[\abs\gm]+\BBK[\abs\gl]\right)h(x)\gr(x)dx<\infty.
\ee
However this type of condition is not satisfactory since it may not hold if $\gm$ and $\gl$ are merely integrable functions since the problem admits always weak solutions. More generally it does not define a  clear class of measures for which we can solve problem $(\ref{I15})$. We introduce new classes of Radon measures whose Green and Poisson potentials belong to a weighted Marcinkiewicz space-type space. Let $\Psi$ be a continuous 
nonincreasing and nonnegative function defined on $[0,\infty)$ and $m$ is a bounded positive Borel measure in $\Gw$ and denote
 \bel{I16'}
M^\Psi_{m}(\Gw):=\left\{f\in \CB(\Gw): \exists \,C>0\;\text { s.t. } \myint{\gl_f(t)}{}dm(x)\leq C\Psi (t)\,,\;\forall t>0\right\}
\ee
where  $\CB(\Gw)$ denotes the space of Borel functions in $\Gw$ and $\gl_f(t)=\{x\in\Gw:\abs{f(x)}>t\}$.
The main results of this article are the two next statements:\smallskip

\noindent{\bf Theorem A }{\it  Let $g$ be an element of  the class $G_{h,\Psi}$ with $\gr h\in L^1(\Gw)$. Then for any $\gm\in \mathfrak M_{\gr}(\Gw)$ and $\gl\in \mathfrak M(\prt\Gw)$ such that $\BBG[\abs\gm]$ and  $\BBK[\abs\gl]$ belong to $M^\Psi_{\gr\,h}(\Gw)$, there exists a solution to problem $(\ref{I15})$. If $r\mapsto g(x,r)$ is nondecreasing for a.e. $x\in\Gw$, this solution is unique.}
\smallskip

Actually we shall introduce a unique formulation for the data $(\gm,\gl)$ as a {\it unique} measure $\gw$ on $\overline\Gw$ which allows to replace $(\ref{I15})$ by $(\ref{I1})$, and a unique assumption on the {\it extended} Green operator $\overline\BBG[\abs\gw]$. We prove in particular the following:\smallskip

\noindent{\bf Theorem B }{\it Assume the assumptions on $h$, $\Psi$ and $g$ of Theorem A are satisfied and $r\mapsto g(x,r)$ is nondecreasing. If $\{(\gw_n\}$ is a sequence of measures in 
$\mathfrak M_{\gr}(\overline\Gw)$ which converges to $\gw\in \mathfrak M_{\gr}(\overline\Gw)$ in the sense that
\bel{I3'}
\myint{\overline\Gw}{}\gz d\gw_n\to\myint{\overline\Gw}{}\gz d\gw
\ee
for any $\gz$ such that $\gr^{-1}\gz\in C(\overline\Gw)$ and if  the $\overline\BBG[\abs{\gw_n}]$ are bounded in $M^\Psi_{\gr\,h}(\Gw)$, then the corresponding solutions $u_{\gw_n}$ of problem $(\ref{I1-n})$
converges to the solution $u_\gw$ of problem $(\ref{I1})$. If $g$ satisfies the $\Gd_2$conditions, the convergence remains valid if only  the $\overline\BBG[\abs{\gw_{s\,n}}]$ are bounded in $M^\Psi_{\gr\,h}(\Gw)$, where $\gw_{s\,n}$ denotes the singular parts of $\gw_n$.
}\smallskip

\section{Linear equations and measures}

Since $\prt\Gw$ is $C^2$, there exists $\gd_{0}>0$ such that, If $x\in\Gw$ is such that $\gr(x)\leq \gd_{0}$, there exists a unique $\gs:=\gs(x)\in \prt\Gw$ such that $\abs{x-\gr(x)}=\gr(x)$. For $\gd>0$ we denote
$$\Gw_{\gd}:=\{x\in\Gw:\gr(x)>\gd\}\,,\;\Gw'_{\gd}:=\{x\in\Gw:\gr(x)<\gd\}\,,\;\Gs_{\gd}:=\{x\in\Gw:\gr(x)=\gd\}\,,\;\Gs:=\Gs_{0}=\prt\Gw.$$
The mapping $x\mapsto (\gr(x),\gs(x))$ is a $C^1$ diffeomorphism from $\overline{\Gw'_{\gd_{0}}}$ onto $[0,\gd_{0}]\ti\Gs$.

 \subsection{Weighted measures on $\overline\Gw$}

We denote by $\mathfrak M(\Gw)$ the set of Radon measures in $\Gw$. If $\ga\in [0,1]$, we denote by $\mathfrak M_{\gr^\ga}(\Gw)$ the subset of $\mathfrak M(\Gw)$ of measures such that
 \bel{M-1}
 \norm{\gm}_{\mathfrak M_{\gr^\ga}}:=\myint{\Gw}{}\gr^\ga d\abs\gm <\infty.
 \ee
 We also set
 \bel{M-2}
C_{\ga}(\overline\Gw):=\{\gz\in C(\Gw):\gr^{-\ga}\gz\in C(\overline\Gw)\} \},
 \ee
 with norm
  \bel{M-3}
\norm\gz_{C_{\ga}}:=\sup_{x\in\overline\Gw}\gr^{-\ga}(x)\abs{\gz(x)}.
 \ee
 Thus, if $\gm\in\mathfrak M_{\gr^\ga}(\Gw)$ and $\gz\in C_{\ga}(\overline\Gw)$, there holds
   \bel{M-4}
\abs{\myint{\Gw}{}\gz d\gm}\leq  \norm{\gm}_{\mathfrak M_{\gr^\ga}}\norm\gz_{C_{\ga}}.
 \ee
 Furthermore, since 
 $$\myint{\Gw_{\gd_0}}{}\gr^\ga d\abs\gm+\sum_{n=1}^\infty\myint{\{2^{-n}\gd_0<\gr\leq 2^{1-n}\gd_0\}}{}\gr^\ga d\abs\gm=\myint{\Gw}{}\gr^\ga d\abs\gm<\infty,
 $$
 there holds
    \bel{M-4-1}
\lim_{\gd\to 0}\myint{\Gw'_\gd}{}\gr^\ga d\abs\gm=0.
 \ee
 We say that a sequence $\{\gm_n\}\subset \mathfrak M_{\gr^\ga}(\Gw)$ converges weakly to $\gm\in \mathfrak M_{\gr^\ga}(\Gw)$ if, for any $\gz\in C_{\ga}(\overline\Gw)$, there holds
   \bel{M-5}
\lim_{n\to\infty}\myint{\Gw}{}\gz d\gm_n=\myint{\Gw}{}\gz d\gm.
 \ee
 However, the  left-hand side expression of $(\ref{M-5})$ may exist but not being a Radon measure in $\Gw$. Therefore 
 we define a more general set  of linear functionals  on $C_{\ga}$
 \bdef {Extended}We denote by $\mathfrak M_{\gr^\ga}(\overline\Gw)$ the set of continuous linear functionals $\gw$ on $C_{\ga}(\overline\Gw)$ such that  there exists a sequence $\{\gm_n\}\subset\mathfrak M_{\gr^\ga}(\Gw)$ which converges weakly to $\gw$.  \es

 The natural norm in $\mathfrak M_{\gr^\ga}(\overline\Gw)$ is
 
   \bel{M-4-1'}
   \norm{\gw}_{\mathfrak M_{\gr^\ga}(\overline\Gw)}=\sup\{\abs{\gw(\gz)}:\gz\in C_{\ga}(\overline\Gw),\norm\gz_{C_\ga}\leq 1\}.
 \ee


 \bprop{property} If $\gw\in \mathfrak M_{\gr^\ga}(\overline\Gw)$, its restriction to $C_c(\Gw)$ is a Radon measure, denoted by $\gm$, which belongs to $\mathfrak M_{\gr^\ga}(\Gw)$.  Furthermore, there exists a Radon measure $\gl$ on $\prt\Gw$ such that 
      \bel{M-6}
\gw (\gz)-\myint{\Gw}{}\gz d\gm=\myint{\prt\Gw}{}\psi\lfloor_{\prt\Gw} d\gl\qquad\forall\gz\in C_{\ga}(\overline\Gw)\text{ and }
\psi=\gr^{-\ga}\gz \in C(\overline\Gw).\ee 
 \es
 \Proof Since $\gw$ is continuous, there exists $C>0$ such that 
       \bel{M-6-1}
\abs{\gw (\gz)}\leq C\norm\gz_{C_\ga}\qquad\forall \gz\in C_{\ga}(\overline\Gw).
 \ee
 This holds in particular if $\gz\in C_{c}(\Gw)$ and proves that the restriction of $\gw$ to $C_{c}(\Gw)$ is a Radon measure that we denote by $\gm$ (as well as the associated Borel measure in $\Gw$) and
 $$\gw (\gz)=\myint{\Gw}{}\gz d\gm\qquad\forall\gz\in C_c(\Gw).
 $$
  Let $\{\gm_n\}\subset \mathfrak M_{\gr^\ga}(\Gw)$ such that 
 $$\lim_{n\to\infty}\myint{\Gw}{}\gz d\gm_n=\gw (\gz)\qquad\forall\gz\in C_\ga(\overline\Gw).
 $$
By the Banach-Steinhaus theorem there exists $C>0$ such that $\norm{\gm_n}_{\mathfrak M_{\gr^\ga}}\leq C$ for all $n\in\BBN$.
 Since for $\gz\in C_c(\Gw)$, 
 $$\gw (\gz)-\myint{\Gw}{}\gz d\gm=\lim_{n\to\infty}\myint{\Gw}{}\gz d(\gm_n-\gm)
 $$
 and
 $$\abs{\myint{\Gw}{}\gz d(\gm_n-\gm)}\leq 2C\norm{\gz}_{C_{\ga}},
 $$
 it follows that $\{\gl_n\}:=\{\gr^\ga(\gm_n-\gm)\}$ is a sequence of Radon measures on $\Gw$, bounded in $\mathfrak M_{\gr^\ga}(\Gw)$ and such that
 $$\lim_{n\to\infty}\myint{\Gw}{}\gz d\gl_n=0\qquad\forall \gz\in C_c(\Gw).
 $$
 Therefore there exists a Radon measure $\gl$ with support in $\prt\Gw$ and a subsequence $\gl_{n_k}$ such that
 $$\lim_{n\to\infty}\myint{\Gw}{}\psi d\gl_{n_k}=\myint{\prt\Gw}{}\psi\lfloor_{\prt\Gw}  d\gl,
 $$
which implies $(\ref{M-6})$.\qeda

\bcor{decomposition} The mapping  $T:\mathfrak M_{\gr^\ga}(\Gw)\ti \mathfrak M(\prt\Gw)\mapsto\mathfrak M_{\gr^\ga}(\overline\Gw) $ 
 defined by
      \bel{M-6-2}
T[\gm,\gl] (\gz)=\myint{\Gw}{}\gz d\gm+\myint{\prt\Gw}{}\psi\lfloor_{\prt\Gw} d\gl\qquad\forall\gz\in C_{\ga}(\overline\Gw)\text{ and }
\psi=\gr^{-\ga}\gz \in C(\overline\Gw).\ee 
is one to one. Furthermore
       \bel{M-6-3}
\max\{\norm\gm_{\mathfrak M_{\gr^\ga}(\Gw)},\norm\gl_{\mathfrak M(\prt\Gw)}\}\leq\norm{T[\gm,\gl]}_{\mathfrak M_{\gr^\ga}(\overline\Gw)}\leq \norm\gm_{\mathfrak M_{\gr^\ga}(\Gw)}+\norm\gl_{\mathfrak M(\prt\Gw)}.
 \ee
\es
\Proof The mapping $T$ is onto from \rprop{property}. The mapping $T$ is one to one since if $T[\gm,\gl]=0$, then $\gm=0$ and $\myint{\prt\Gw}{}\psi\lfloor_{\prt\Gw} d\gl=0$ for any $\psi\in C(\overline\Gw)$. This implies $\gl=0$. The right-hand side inequality $(\ref{M-6-3})$ is clear since $\sup\abs{\psi\lfloor_{\prt\Gw}}\leq \norm\gz_{C_\ga}$. Because of $(\ref{M-4-1})$
$$\myint{\Gw}{}\gr^\ga d\abs{\gm}=\sup\left\{\myint{\Gw}{}\gz d\gm:\gz\in C_c(\Gw), \norm\gz_{C_\ga}\leq 1\right\}
$$
This implies 
$$\norm\gm_{\mathfrak M_{\gr^\ga}(\Gw)}\leq\norm{T[\gm,\gl]}_{\mathfrak M_{\gr^\ga}(\overline\Gw)}
$$
If $\phi\in C(\prt\Gw)$ is such that $\abs\phi\leq 1$ and $\Phi$ is its harmonic lifting in $\Gw$, the function $\gz=\gr^\ga\Phi$ belongs to $C_\ga(\overline\Gw)$ and satisfies $\norm{\gz}_{C^\ga}\leq 1$. Let $\{\eta_n\}\subset C^\infty(\BBR^N)$ such that $0\leq\eta_n\leq 1$, $\eta_n(x)=0$ if $\gr(x)\geq 2/n$, $\eta_n(x)=1$ if $\gr(x)\leq 1/n$. Then $\gz_n=\eta_n\gr^\ga\Phi$ belongs also to $C_\ga(\overline\Gw)$ and $\norm{\gz_n}_{C^\ga}\leq 1$. Since
$$T[\gm,\gl] (\gz_n)=\myint{\Gw}{}\gz_n d\gm+\myint{\prt\Gw}{}\phi d\gl
$$
and $\myint{\Gw}{}\gz_n d\gm\to 0$ as $n\to\infty$, we derive
$$\norm{T[\gm,\gl]}_{\mathfrak M_{\gr^\ga}(\overline\Gw)}\geq \myint{\prt\Gw}{}\phi d\gl.
$$
This ends to proof.\qeda\medskip

\noindent\Remark If $\gl$ is a Radon measure on $\prt\Gw$ and we can define its $\gd^\ga$-lifting $\Gl_{\gd^\ga}[\gl]\in \mathfrak M(\Gw)$ by
$$\myint{\Gw}{}\gz d\gl_{\gd^\ga}=\gd^{-\ga}\myint{\Gw}{}\gz(\gd,\gs)d\gl(\gs).
$$
Clearly $\gl_{\gd^\ga}\in \mathfrak M_{\gr^\ga}(\Gw)$ and if $\gz\in C_{\ga}(\overline\Gw)$ and 
 $\ell_{\ga}(\gz)=-\lim_{\gr\to 0}\gr^{-\ga}\gz$, then $\ell_{\ga}(\gz)\in C(\prt\Gw)$, there holds
       \bel{M-8}
\lim_{\gd\to 0}\myint{\Gw}{}\gz d\gl_{\gd^\ga}=\myint{\Gs}{}\ell_{\ga}(\gz) d\gl.
 \ee
In the particular case where $\ga=1$ $\ell_{\ga}(\gz)=\gz_{{\bf n}}:=\lim_{\gr\to 0}\gr^{-1}\gz$, and
      \bel{M-9}
\lim_{\gd\to 0}\myint{\Gw}{}\gz d\gl_{\gd}=-\myint{\Gs}{}\gz_{{\bf n}}d\gl.
 \ee
\subsection{The linear operator}
Let $x=(x_1,...,x_N)$ the coordinates in $\BBR^N$ and $\Gw$ a bounded domain in $\BBR^N$. We consider the operator $L$
 in divergence form defined by
 \bel{2-2-1}
Lu:=-\sum_{i,j=1}^N\myfrac{\prt}{\prt x_i}\left(a_{ij}\myfrac{\prt u}{\prt x_j}\right)+\sum_{i=1}^Nb_{i}\myfrac{\prt u}{\prt x_i}-
\sum_{i=1}^N\myfrac{\prt}{\prt x_i}\left(c_{i}u\right) +du
\ee
where the $a_{ij}$, $b_i$ and $c_i$  are Lipschitz continuous and $d$ is bounded and measurable in $\Gw$. We assume that the ellipticity condition
 \bel{2-2-2}
\sum_{i,j=1}^Na_{ij}(x)\xi_i\xi_j\geq a\sum_{i1}^N\xi_i^2\qquad\forall\xi\in\BBR^N
\ee
holds for almost $x$ in $\Gw$, for some $a>0$. We also assume the positivity condition
 \bel{2-2-3}
\int_\Gw\left(dv+\myfrac{1}{2}\sum_{i=1}^N(b_i+c_i)\myfrac{\prt v}{\prt x_i}\right)dx\geq 0\qquad\forall v\in C^1_c(\Gw),\,v\geq 0
\ee
Under these assumptions, the bilinear form
 \bel{2-2-4}
(u,v)\mapsto A_L(u,v)=\int_\Gw\left(\sum_{i,j=1}^Na_{ij}\myfrac{\prt u}{\prt x_j}\myfrac{\prt v}{\prt x_i}
+\sum_{i=1}^N\left(b_i\myfrac{\prt u}{\prt x_i}v+c_i\myfrac{\prt v}{\prt x_i}u\right)+duv\right)dx
\ee
is continuous and coercive on $W^{1,2}(\Gw)$. We define the adjoint operator  $L^*$ by
 \bel{2-2-5}
L^* u:=-\sum_{i,j=1}^N\myfrac{\prt}{\prt x_j}\left(a_{ij}\myfrac{\prt u}{\prt x_i}\right)+\sum_{i=1}^Nc_{i}\myfrac{\prt u}{\prt x_i}-
\sum_{i=1}^N\myfrac{\prt}{\prt x_i}\left(b_{i}u\right) +du
\ee

We denote by $G=G_L$ and $K=K_L$ the Green and Poisson kernels corresponding to the operator $L$ in $\Gw$. We recall the following equivalence statement \cite{Pin}, \cite{An} 
\bprop{equiv} Assume $\Gw$ has a $C^2$ boundary and $(\ref{2-2-3})$ holds. Then there exists a positive constant $C$ such that
 \bel{2-2-6}
CG_{-\Gd}\leq G\leq C^{-1}G_{-\Gd}\qquad\text{in }\Gw\ti\Gw\setminus D_\Gw
\ee
where $D_\Gw={x\in\Gw\ti\Gw:x:\neq y}$ and
 \bel{2-2-7}
CK_{-\Gd}\leq K\leq C^{-1}K_{-\Gd}\qquad\text{in }\Gw\ti\prt\Gw.
\ee
\es
 \subsection{Linear equation with measure data}
 If $m\in\mathfrak M_+(\Gw)$ is a bounded  Borel measure and $\Psi:[0,\infty)\mapsto [0,\infty)$ is continuous and nonincreasing, we define the subset $M_{m}^\Psi(\Gw)$ of the set $\CB(\Gw)$ of Borel mesurable functions by
 \bel{2-3-1}
M^\Psi_{m}(\Gw):=\left\{f\in \CB(\Gw):\exists C>0\;\text{ s.t. }\myint{\gl_{f}(t)}{} dm(x)\leq C\Psi (t)\,,\;\forall t>0\right\}
\ee
where 
 \bel{rep}
 \gl_{f}(t)=\{x\in\Gw:\abs{f(x)}>t\}.
\ee
Notice that $\Psi (t)\leq m(\Gw)$ for $t\geq 0$. Denote
 \bel{rep'}
\bar \gl_{f}(t)=\{x\in\Gw:\abs{f(x)}\geq t\}.
\ee
Since $\Psi$ is continuous, $(\ref{2-3-1})$ implies 
$$\myint{\bar\gl_{f}(t)}{} dm(x)\leq C\Psi (t)\,,\;\forall t>0. 
$$
If we modify $\Psi$ in order to impose $\Psi (0)= m(\Gw)$, $(\ref{2-3-1})$ is equivalent to
 \bel{2-3-1'}
M^\Psi_{m}(\Gw):=\left\{f\in \CB(\Gw):\exists C>0\;\text{ s.t. }\myint{\bar\gl_{f}(t)}{} dm(x)\leq C\Psi (t)\,,\;\forall t\geq 0\right\}
\ee

We denote by $C^\Psi_{m}(f)$ the smallest constant $C$ such that $(\ref{2-3-1})$ holds. If $t\mapsto \Psi (t)/\Psi(2t)$ remains bounded on $[0,\infty)$, $M^\Psi_{m}(\Gw)$ is a vector space $f\mapsto C^\Psi_{m}(f)$ is a quasi-norm on 
the quotient space $M^\Psi_{m}(\Gw)/\CR$ where $\CR$ is the equivalence relation $f_1\CR f_2\Longleftrightarrow 
f_1-f_2=0\;m\text {-a.e. in }\; \Gw$. In general $M^\Psi_{m}(\Gw)$ is not a vector space

When $\Psi (t)=t^{-p}$ with $p\geq 1$ and $m (x)=\gr(x)^\ga$, with $\ga\in [0,1]$, we denote by $M^p_{\gr^\ga}(\Gw)$ the corresponding Marcinkiewicz space.  The following results proved in \cite{BVVi} with $L=-\Gd$ are valid for a general operator $L$

 \bprop{reg} Let $\ga\in [0,1]$, $N\geq 2$. If $\gm\in \mathfrak M_{\gr^\ga}(\overline\Gw)$ and $N+\ga-2>0$,  
   \bel{2-3-2}
\norm{\BBG[\gm]}_{M^{(N+\ga)/(N+\ga-2)}_{\gr^\ga}}\leq C\norm{\gm}_{\mathfrak M_{\gr^\ga}},
\ee
   \bel{2-3-3}
\norm{\nabla\BBG[\gm]}_{M^{(N+\ga)/(N+\ga-1)}_{\gr^\ga}}\leq C\norm{\gm}_{\mathfrak M_{\gr^\ga}}.
\ee
Furthermore, for any $\gamma\in [0,1]$ and $\gl\in \mathfrak M(\prt\Gw)$,
   \bel{2-3-4}
\norm{\BBK[\gl]}_{M^{(N+\gamma)/(N-1)}_{\gr^\gamma}}\leq C\norm{\gl}_{\mathfrak M}.
\ee
\es
We recall the following result proved in \cite[Th 2.9]{Ve1}
\bth {Lin} Let $\ga\in [0,1]$. For every $\gm\in \mathfrak M_{\gr^\ga}(\Gw)$ and $\gl\in \mathfrak M(\prt\Gw)$, there exists a unique $u:=u_{\gm,\gl}\in L^1(\Gw)$ satisfying
   \bel{2-3-5}\BA {ll}
   -Lu=\gm\qquad&\text{in }\Gw\\
   \phantom{-L}
u=\gl\qquad&\text{in }\prt\Gw,
\EA\ee
in the following weak sense
   \bel{2-3-6}\BA {l}
-\myint{\Gw}{}uL^*\gz dx=\myint{\Gw}{}\gz d\gm-\myint{\prt\Gw}{}\gz_{\bf n} d\gl\qquad\forall \gz\in C_{c^{1,L}}(\overline\Gw).
\EA\ee
Furthermore, if $\{(\gm_{n,\gl_{n}})\}$ is bounded in $\mathfrak M_{\gr^\ga}(\Gw)\ti\mathfrak M(\prt\Gw)$ and converges weakly with respect to $C_\ga(\overline\Gw)\ti C(\prt\Gw)$ to $(\gm,\gl)\in \mathfrak M_{\gr^\ga}(\Gw)\ti\mathfrak M(\prt\Gw)$, then $u_{{\gm_{n},\gl_{n}}}$ converges to $u_{\gm,\gl}$.
\es

\noindent\Remark If we define the measure $\gw\in \mathfrak M_{\gr^\ga}(\overline\Gw)$ by $\gw=T[\gm,\gl]$ (see $(\ref{M-6-2})$), then it can also be expressed by
   \bel{2-3-7}\BA {l}
\myint{\overline\Gw}{}\gz d\gw:=\myint{\Gw}{}\gz d\gm-\myint{\prt\Gw}{}\gz_{\bf n} d\gl\quad\forall\gz\in C_{1}(\overline\Gw),
\EA\ee
since $\gz\in C_{1}(\overline\Gw)$ implies that $\gz_{\bf n}$ exists on $\prt\Gw$ and is continuous. We define the {\it global } Green operator on $\overline\Gw$ by
   \bel{2-3-8}\BA {l}
\overline\BBG[\gw]:=\BBG[\gm])+\BBP_{L}[\gl].
\EA\ee
and $(\ref{2-3-5})$ is replaced by the unique equation
   \bel{2-3-5'}
 -Lu=\gw\qquad\text{in }\overline\Gw.
    \ee
Then $(\ref{2-3-2})$-$(\ref{2-3-4})$ with $\ga=1$ are equivalent to
   \bel{2-3-9}
\norm{\overline\BBG[\gw]}_{M^{(N+1)/(N-1)}_{\gr}}\leq C\norm{\gw}_{\mathfrak M_{\gr}}.
\ee
Furthermore, we say that $u\in L^1(\Gw)$ is a subsolution of $(\ref{2-3-5'})$ in $\overline\Gw$, if 
\bel{2-3-10}
-\myint{\Gw}{}uL^*\gz dx\leq \myint{\overline\Gw}{}\gz d\gw:=\myint{\Gw}{}\gz d\gm-\myint{\prt\Gw}{}\gz_{\bf n} d\gl\qquad\forall \gz\in C^{1,L^*}_c(\overline\Gw)\,,\;\gz\geq 0.
\ee
Comparison principle applies, thus $u\leq \overline \BBG[\gw]$. A supersolution is defined similarly.\medskip

\noindent\Remark If $\gw=T[\gm,\gl]\in \mathfrak M^+_{\ga}(\overline\Gw)$ its Lebesgue decomposition is 
$\gw_r+\gw_s=T[\gm_r,\gl_r]+T[\gm_s,\gl_s]$ where $\gm_r$ and $\gl_r$ are the absolutely continuous part with respect to the Hausdorff measures $d\CH^N$ and $d\CH^{N-1}$ and $\gm_s$ and $\gl_s$ the respective singular parts. Similarly if $\gw=T[\gm,\gl]$, then $\gw=\gw^+-\gw^-$ where $\gw^+=T[\gm^+,\gl^+]$ and $\gw^-=T[\gm^-,\gl^-]$.
\subsection{Regularity results}

We define the class of measures $B^p_{h}(\overline\Gw)$ by
\bel{2-4-1}
B^\Psi_{h}(\overline\Gw):=\{\gw\in \mathfrak M_{\gr}(\overline\Gw):\overline\BBG[\abs{\gw}]\in M^\Psi_{\gr h}\Gw)\}.
\ee
By \rprop {equiv}, this class remains unchanged if we replace $-\Gd$ by $L$ and the Green operator for $L$ by the one of $-\Gd$.  
If $\Psi(t)=t^{-p}$ and $h=1$, the corresponding class of measures is larger that the usual 
    \bel{2-4-2}
\tilde B^p(\overline\Gw):=\{\gw\in \mathfrak M_{\gr}(\overline\Gw):\overline\BBG[\abs{\gw}]\in L^p_{\gr}(\Gw)\}
\ee
which corresponds to negative Besov spaces: if $\gw=T[\gm,\gl]$, then the regularity results for harmonic functions  \cite {MV1} and solution of Laplace equation \cite {AH} yields to
    \bel{2-4-1'}
\tilde B^p(\overline\Gw)=B^{-\frac{2}{p},p}(\Gw).
\ee 

\noindent {Example 1} If $h(x)=(\gr(x))^\gb$, with $\gb>-2$. Then $\gw=T[0,\gl]\in B^p_{\gr^\gb}(\overline\Gw)$ if and only if
$\overline\BBG[|\gw|]\in M_{\gr^{\gb+1}}(\Gw)$. This means that  $\gl\in B^{-s,p}_\infty(\prt\Gw)$ with $s=(\gb+2)/p$ (see \cite{Tri} for the definition of $B^{\ga,p}_q$.\medskip


\section{The main results}

\bdef {class}We say that  a Caratheodory function $g:\Gw\ti\BBR$ belongs to the class $G_{h,\Psi}$ if there exist a nonnegative function $h\in L^1_\gr(\Gw)$, a continuous nondecreasing function $\tilde g$ defined on $\BBR_+$ and vanishing at $r=0$  such that 
$0\leq g(x,r)\rm {sign\,}r\leq h(x)\tilde g(\abs r)$ in $\Gw\ti\BBR$ and a continuous nonincreasing function $\Psi:[0,\infty)\mapsto [0,\infty)$ with the property that
    \bel{3-1}
-\myint{1}{\infty}\tilde g(s)d\Psi(s)<\infty.
\ee
\es 
\blemma{equiv} Let $\gm$ be a nonnegative measure in $\mathfrak M(\Gw)$ and $g:\Gw\ti\BBR\mapsto \BBR$ a Caratheodory function such that $0\leq g(x,r)\rm {sign\,}r\leq h(x)\tilde g(\abs r)$ where $h\in L^1_\gr(\Gw)$ and $\tilde g$ is  a continuous and nondecreasing function $\tilde g$ defined on $\BBR_+$ and vanishing at $r=0$. Then \smallskip

\noindent (i) If $g\in G_{h,\Psi}$ and $\gm\in B^\Psi_{h}(\overline\Gw)$, then $\tilde g\circ \overline\BBG[\gm]\in L^1_{\gr h}(\Gw)$.\smallskip

\noindent (ii) if $\tilde g\circ \overline\BBG[\gm]\in L^1_{\gr h}(\Gw)$ and , then $\gm\in B^\Psi_{h}(\overline\Gw)$ and $g\in G_{h,\Psi}$ with $\Psi(s)=\gth_{\gl_{\overline\BBG[\gm]}}(s)$, where $\gl_{\overline\BBG[\gm]}(s)$ is defined by $(\ref{rep})$ with $f$ replaced by  $\overline\BBG[\gm]$ and $\gth_{\gl_{\overline\BBG[\gm]}}(s)=\myint{\gl_{\overline\BBG[\gm]}(s)}{}d(\gr h)$.
\es
\Proof This due to the fact that
 \bel{eq}
   \myint{\Gw}{}\tilde g(\overline\BBG[\gm])\gr h dx= -\myint{0}{\infty}\tilde g(s)d\gth_{\gl_{\overline\BBG[\gm]}}(s).
  \ee
  Therefore, if $\gth_{\gl_{\overline\BBG[\gm]}}(s)\leq \Psi(s)$, it proves (i).
Conversely, if $\Psi(s)=\gth_{\gl_{\overline\BBG[\gm]}}(s)$, then  $\gm\in B^\Psi_{h}(\overline\Gw)$ and   $g\in G_{h,\Psi}$.
 \qeda\medskip
 
The following existence result extends to one in \cite{Ve1}
\bth{Exist} Assume $g$ belongs to the class $G_{h,\Psi}$. Then for any  $\gw\in B^\Psi_{h}(\overline\Gw)$ there exists a function $u\in L^1(\Gw)$ such that $g\circ u\in L^1(\Gw)$ satisfying
    \bel{3-2}
\myint{\Gw}{}\left(-uL^*\gz+g\circ u\,\gz\right)dx=\myint{\overline\Gw}{}\gz d\gw\qquad\forall\gz\in C_c^{1,L^*}(\overline\Gw).
\ee
Furthermore $u$ is unique if $r\mapsto g(x,r)$ is nondecreasing for a.e. $x\in\Gw$.
\es
\Proof It is essentially \cite[Theorem 3.7]{Ve1}. Since $0\leq g(x,r)\rm {sign\,}r\leq h(x)\tilde g(\abs r)$, we define the following truncation  $g_{k}(.,r)$ for any $k>0$.
\begin {equation}\label {trunc}
g_{k}(x,r)=g(x,r)\chi_{{\Gth_k}}
\end{equation}
where $\Gth_k=\{x\in\Gw:h(x)\leq k\}$. Then $0\leq g(x,r)\rm {sign\,}r\leq k\tilde g(\abs r)$ and there exists a solution $u_k$ to
\bel{X1}\BA {ll}
-Lu_k+g_k\circ u_k=\gw\qquad&\text{in }\overline\Gw\EA.\ee
Actually, in \cite[Theorem 3.7]{Ve1} the proof is done with $\gm\in \mathfrak M_{\gr^\ga}(\Gw)$ for any $\ga\in [0,1]$, but due to our definition of measures in $\mathfrak M_{\gr^\ga}(\overline\Gw)$, it is also valid in this case.
\smallskip



\noindent {\it Step 2: Convergence when $k\to\infty$.}  By Brezis'estimates (see e.g. \cite[Th 2.4]{Ve1}), for any 
 $\gz\in C_{c}^{1,L}(\overline\Gw)$, $\gz\geq 0$, one has
\begin {eqnarray}\label {SDLM2b}
\int_{\Gw}\left(- \abs {u_{k}}L^*\gz +{\rm 
sign(u_{k})}g_k(x,u_{k})\gz\right)dx\leq\int_{\overline\Gw}\gz d\abs{\gw}.
\end {eqnarray} 
and
\bel{SDLM3b}
\norm {u_{k}}_{L^1}+\norm {\gr g_k(.,u_{k})}_{L^1_\gr}
\leq C_{1}\norm {\gw}_{\mathfrak M_\gr}.
\ee
Furthermore, by estimates of \rprop{reg} and since $\abs{u_{k}}\leq \overline\BBG[\abs{\gw}]$, there holds, 
\bel{SDLM4}
\norm{u_{k}}_{M^{(N+1)/N}_{\gr}}+\norm{\nabla u_{k}}_{M^{(N+1)/N}_{\gr}}
\leq C\norm {\gw}_{\mathfrak M_{\gr}}.
\ee
Since the right-hand side of $(\ref{SDLM4})$ is bounded independently of  $k$ fixed, there exist a subsequence $\{u_{k_j}\}$ and a function $u\in W^{1,q}_{loc}(\Gw)$, 
for any $1\leq q<(N+1)/N$, such that $u_{k_j}\to u$ a.e. in $\Gw$ - and thus $g_{k_j}\circ u_{k_j}\to g\circ u$ a.e. - and 
weakly in  $W^{1,q}_{loc}(\Gw)$ when $k_j\to\infty$.
Let $R>0$ and $E\subset\Gw$ be a Borel subset, then
\bel{SDLM4b'}\BA {l}
\myint{E}{}\abs {g_{k_j}\circ u_{k_j}}\rho dx
\leq \myint{E\cap\{\abs{u_{k_j}}\leq R\}}{}\tilde g(\abs{u_{k_j}})\rho hdx
+\myint{E\cap\{\abs{u_{k_j}}> R\}}{}\tilde g(\abs{u_{k_j}})\rho hdx\\[4mm]
\phantom{\myint{E}{}\abs {g_{k_j}\circ u_{k_j}}\rho dx}
\leq \tilde g(R)\myint{E}{}\gr hdx
-\myint{R}{\ity}\tilde g( s)d\gth_{u_{k_j}}(s),
\EA
\ee
where, we recall it, 
\begin {eqnarray*}
\gth_{u_{k_j}}(s):=\myint{\gl_{u_{k_j}}(s)}{}d(\gr h).
\end {eqnarray*}
Since $\abs {u_{k_j}}\leq \overline\BBG[\abs{\gw}]$, $\gth_{u_{k_j}}(s)\leq \gth_{\overline\BBG[\abs{\gw}]}(s)$. By assumption, 
$$\gth_{\overline\BBG[\abs{\gw_n}]}(s)\leq C\Psi (s)\qquad\forall s>0,
$$
with 
$$C=C^\Psi_{\gr h}(\overline\BBG[\abs{\gw}]).
$$
Furthermore, by a standard integration by parts in Stieltjes integrals and for a.e. $R$,
\bel{SDLM4b}\BA {l}
-\myint{R}{\ity}\tilde g( s)d\gth_{u_{k_j}}(s)=\tilde g(R)\gth_{u_{k_j}}(R)+
\myint{R}{\ity}\gth_{u_{k_j}}(s)d\tilde g( s))\\
[4mm]\phantom{-\myint{R}{\ity}\tilde g( s)d\gth_{u_{k_j}}(s)}
\leq \tilde g(R)\gth_{u_{k_j}}(R)+
C\myint{R}{\ity} \Psi(s)d\tilde g( s)\\
[4mm]\phantom{-\myint{R}{\ity}\tilde g( s)d\gth_{u_{k_j}}(s)}
\leq\tilde g(R)\gth_{u_{k_j}}(R)-C\tilde g(R)\Psi(R) -C\myint{R}{\ity}\tilde g( s)d\Psi(s)\\
[4mm]\phantom{-\myint{R}{\ity}\tilde g( s)d\gth_{u_{k_j}}(s)}
\leq -C\myint{R}{\ity}\tilde g( s)d\Psi(s).
\EA\ee
Since condition $(\ref{3-1})$ holds, it follows
\begin {equation}\label {marsest2b}
\lim_{R\to\infty}\myint{R}{\ity}\tilde g( s)d\Psi(s)=0.
\end {equation}
Given $\epsilon>0$, we first choose $R>0$ such that 
$$-C\myint{R}{\ity}\tilde g( s)d\Psi(s)\leq \epsilon/2.
$$
Then we put $\gd=\epsilon/(2(1+\tilde g(R))$ and derive
$$\int_{E}\rho  dx\leq \gd\Longrightarrow 
\int_{E}\abs {g_{k_j}(u_{k_j})}\rho h dx\leq \epsilon.
$$
Therefore $\{g_{k_j}\circ u_{k_j}\}$ 
is uniformly integrable in $L^1_\gr(\Gw)$. It follows by Vitali's convergence theorem
\begin {equation}\label {lim}
\lim_{k\to\ity}g_{k_j}\circ u_{k_j}=g\circ u\qquad\text{in }L^1_\gr(\Gw).
\end {equation}
Let $\gz\in C_{c}^{1,L}(\overline\Gw)$. If we let $k_j\to\infty$ in the equality
\bel{X1'}
\int_{\Gw}\left(- u_{k_j}L^*\gz +g_{k_j}\circ u_{k_j}\gz\right)dx=\int_{\overline\Gw}\gz d\gw,
\ee
we derive 
\bel{X2}
\int_{\Gw}\left(- uL^*\gz +g\circ u\gz\right)dx=\int_{\overline\Gw}\gz d\gw.
\ee
Uniqueness follows classicaly if $g(x,.)$ is nonndecreasing.
\qeda \medskip

The following extension of the previous result is an adaptation of \cite[Th. 3.20]{Ve1}

\bth {delta2}Assume $g$ belongs to the class $G_{h,\Psi}$  and satisfies the following $\Gd_2$-condition
\bel {Delta2}
\abs{g(x,r+r')}\leq \gth\left(\abs{g(x,r)}+\abs{g(x,r')}\right)+\ell(x)\qquad\forall x\in\Gw,\,\forall (r,r')\in\BBR\ti\BBR,
\ee
for some nonnegative $\ell \in L^1_\gr(\Gw)$. Suppose also that $r\mapsto g(x,r)$ is nondeacreasing.  If $\gw\in \mathfrak M_\gr(\overline\Gw)$ has Lebesgue decomposition $\gw=\gw_r+\gw_s$ with  regular part with respect to the Lebesgues measures $\gw_r$ and singular part  $\gw_s$, and if $\gw_s$ belongs to $B^\Psi_{h}(\overline\Gw)$, then there exists a unique solution $u$ to $(\ref{3-2})$.
\es
\Proof If $g$ satisfies $(\ref{Delta2})$, $g_k$ defined by $(\ref{trunc})$ shares the same property with the same $\ell$. Therefore, by \cite[Th 3.12]{Ve1}, there exists a solution $u_k$ to $(\ref{X1})$. Actually, in this result it is only assume that 
$\ell$ in $(\ref{Delta2})$ is a constant, but the proof is valid if it is a nonnegative  function in $L^1_\gr(\Gw)$. Let  $v_k$ and $v'_k$ be weak solutions in $\overline\Gw$ of 
$-Lv_k+g_k\circ v_k=\gw^+_r$ and $-Lv_k'-g_k\circ (-v_k')=\gw^-_r$ respectively. Set $w_k=v_k+\overline\BBG(\gw^+_s)$ and $w_k'=v_k'+\overline\BBG(\gw^-_s)$. Then $-Lw_k+g_k\circ w_k\geq \gw^+$ and $-Lw_k'-g_k\circ (-w_k')\geq \gw^-$ in $\overline\Gw$. By monotonicity
$-w'_k\leq u_k\leq w_k $, thus $g_k(-w'_k)\leq g_k(u_k)\leq g_k(w_k)$. The estimates $(\ref{SDLM3b})$ and $(\ref{SDLM4})$ are satisfied, therefore there exist a function $u\in L^1(\Gw)$ and a subsequence $u_{k_j}$ which converges to $u$ a.e. in 
$\Gw$. Furthermore
\bel{Delta3}\BA {l} g_k(x,u_k)\leq \gth\left(g_k(x,v_k)+g_k(x,\overline\BBG(\gw^+_s)\right)+\ell\\
\phantom{g_k(x,u_k)}\leq  \gth\left(g_k(x,v_k)+g(x,\overline\BBG(\gw^+_s)\right)+\ell
\EA\ee
Since the sequence $\{\abs{g_k}\}$ increases, $\{v_k\}$ and $\{v'_k\}$ decrease. Therefore $v_k\downarrow v$  and $v'_k\downarrow v'$ which satisfy $-Lv+g\circ v=\gw^+_r$ and $-Lv'-g_k\circ (-v')=\gw^-_r$ respectively in $\overline\Gw$. Therefore $g_k\circ v_k\to g\circ v$
and $g_k\circ v'_k\to -g\circ (-v')$ in $L^1_\gr(\Gw)$ respectively. Since 
$$g_k\circ \overline\BBG(\gw^+_s)\leq g\circ \overline\BBG(\gw^+_s)$$
and $\gw_s\in B^\Psi_{h}(\overline\Gw)$, $g\circ \overline\BBG(\gw^+_s$ by \rlemma{equiv}, the right-hand side term of inequality $(\ref{Delta3})$ is uniformly integrable in $L^1_\gr(\Gw)$. Similarly 
\bel{Delta4}g_k(x,u_k)\geq \gth\left(g_k(x,-v'_k)+g(x,-\overline\BBG(\gw^-_s)\right)-\ell
\ee
and the right-hand side of $(\ref{Delta4})$ is also  uniformly integrable in $L^1_\gr(\Gw)$. We conclude as in \rth{Exist}
.\qeda
\section{Stability}
\blemma{closure} Let $\{\gw_n\}\subset B^\Psi_{h}(\overline\Gw)$ be a sequence of measures such that $C^\Psi_{\gr}(\overline\BBG[\abs{\gw_n}])$ is bounded independently of $n$. Then  $\{\gw_n\}$ remains bounded in $\mathfrak M_\gr(\overline\Gw)$. If  $\gw_n\to \gw$ weakly in $\mathfrak M_\gr(\overline\Gw)$, then $\gw\in  B^\Psi_{h}(\overline\Gw)$.
\es
\Proof Since $C^\Psi_{\gr}(\overline\BBG[\abs{\gw_n}])$ is uniformly bounded, the sequence $\{g\circ\overline\BBG[\abs{\gw_n}])\}$ is bounded in $L^1_\gr(\Gw)$ by \rlemma{equiv}. Since $\gw_n\to \gw$ weakly in $\mathfrak M_\gr(\overline\Gw)$, $\overline\BBG[\gw_n]\to\overline\BBG[\gw]$ in $L^1_\gr(\Gw)$ and, up to a subsequence, a.e. in $\Gw$. Therefore, and up to sets of zero Lebesgue measure,

\bel{conv0}
\gl_{\overline\BBG[\gw]}(t)\subset\bigcap_{n\geq 0}\left(\bigcup_{p\geq n}\gl_{\overline\BBG[\gw_p]}(t)\right)\subset
\bigcap_{n\geq 0}\left(\bigcup_{p\geq n}\overline\gl_{\overline\BBG[\gw_p]}(t)\right)
\subset\overline\gl_{\overline\BBG[\gw]}(t).
\ee
Therefore 
\bel{conv0-1}
\limsup_{n\to\infty}\gth_{\gl_{\overline\BBG[\gw_n]}(t)}\leq \gth_{\overline\gl_{\overline\BBG[\gw]}(t)}.
\ee
Conversely, for any $x\in \gl_{\overline\BBG[\gw]}(t)$, i.e. such that $\overline\BBG[\gw](x)>t$, there exists $n_x$ such that $x\in \gl_{\overline\BBG[\gw_n]}(t)$ if $n\geq n_x$. This implies
$$\lim_{n\to\infty}\chi_{_{ \gl_{\overline\BBG[\gw_n]}(t)}}\chi_{_{ \gl_{\overline\BBG[\gw]}(t)}}=\chi_{_{ \gl_{\overline\BBG[\gw]}(t)}},
$$
and
\bel{conv0-2}
\liminf_{n\to\infty}\gth_{ \gl_{\overline\BBG[\gw_n]}(t)}\geq \gth_{ \gl_{\overline\BBG[\gw]}(t)}.
\ee
Since $\gth_{ \gl_{\overline\BBG[\gw_n]}(t)}\leq C^\Psi_{\gr}(\overline\BBG[\abs{\gw_n}])\Psi(t)$ and the $C^\Psi_{\gr}(\overline\BBG[\abs{\gw_n}])$ are bounded, it follows that $\gw$ belongs to $B^\Psi_{h}(\overline\Gw)$.\qeda

\bth{Stab1} Assume $g$ belongs to the class $G_{h,\Psi}$ and  $r\mapsto g(x,r)$ is nondecreasing for a.e. $x\in\Gw$.  Let $\{\gw_n\}\subset B^\Psi_{h}(\overline\Gw)$ be a sequence of measures such that $C^\Psi_{\gr}(\BBG[\abs{\gw_n}])$
is bounded independently of $n$ which converges to $\gw$ weakly with respect to $C_1(\overline\Gw)$. Then the solution $u_n$ of 
\bel{conv1}-Lu_n+g\circ u_n=\gw_n\qquad\text{in }\overline\Gw
\ee
converges to the solution $u$ of
\bel{conv2}-Lu+g\circ u=\gw\qquad\text{in }\overline\Gw
\ee
\es
\Proof Since $u_n$ satisfies the Brezis estimates $(\ref{SDLM3b})$ and $(\ref{SDLM4})$, there exists a subsequence $\{u_{n_j}\}$ 
and $u\in L^1(\Gw)$ such that $u_{n_j}\to u$ a.e. in $\Gw$ and in $L^1(\Gw)$. As in the proof of \rth{Exist}, the problem is to prove the convergence of the $g\circ u_{n_j}$ in $L^1_\gr(\Gw)$. But this is a clearly obtained by the uniform integrability, as in the proof of \rth{Exist}-Step 2, using the fact that, in $(\ref{SDLM4b})$, the $\gth_{u_{n_j}}$ are bounded by $sup_{n}C_{\gr h}^\Psi (\overline\BBG[\gw_n])\Psi$. \qeda

\bth{Stab2}Assume $g$ belongs to the class $G_{h,\Psi}$, satisfies the  $\Gd_2$-condition $(\ref{Delta2})$ and $r\mapsto g(x,r)$ is nondeacreasing.  Let $\{\gw_n\}\subset \mathfrak M_\gr(\overline\Gw)$ has Lebesgue decomposition $\gw_n=\gw_{n\,r}+\gw_{n\,s}$  if $\{\gw_{n\,s}\}\subset B^\Psi_{h}(\overline\Gw)$ are such that the $C_{\gr h}^\Psi (\overline\BBG[\gw_{n\,s}])$ are uniformly bounded, then the solutions $u_{n}$ of $(\ref{conv1})$ converges in $L^1(\Gw)$ to the solution $u$ of $(\ref{conv2})$.
\es
\Proof The argument follows the one of \rth{delta2}.  Let  $v_n$ and $v'_n$ be weak solutions  in $\overline\Gw$ of 
$-Lv_n+g\circ v_n=\gw^+_{n\,r}$ and $-Lv_n'-g\circ (-v_n')=\gw^-_{n\,r}$ respectively. Set $w_n=v_n+\overline\BBG(\gw^+_{n\,s})$ and $w_k'=v_k'+\overline\BBG(\gw^-_{n\,s})$. Then $-Lw_n+g\circ w_n\geq \gw_n^+$ and $-Lw_n'-g\circ (-w_n')\geq \gw_n^-$. By monotonicity
$-w'_n\leq u_n\leq w_n $, thus $g(-w'_n)\leq g(u_n)\leq g(w_n)$. The estimates $(\ref{SDLM3b})$ and $(\ref{SDLM4})$ are satisfied therefore there exist a function $u\in L^1(\Gw)$ and a subsequence $u_{n_j}$ which converges to $u$ a.e. in 
$\Gw$ and in $L^1(\Gw)$. Furthermore
\bel{Delta3'}\BA {l} g(x,u_n)\leq \gth\left(g(x,v_n)+g(x,\overline\BBG(\gw^+_{n\,s})\right)+\ell\\
\phantom{g(x,u_n)}\leq  \gth\left(g(x,v_n)+g(x,\overline\BBG(\gw^+_{n\,s})\right)+\ell.
\EA\ee
Classicaly $v_n\to v$ $v'_n\to v'$ in $L^1(\Gw)$ which satisfy $-Lv+g\circ v=\gw^+_r$ and $-Lv'-g_k\circ (-v')=\gw^-_r$ respectively. Therefore $g\circ v_n\to g\circ v$
and $g\circ v'\to -g\circ (-v')$ in $L^1_\gr(\Gw)$ respectively. Since $C_{\gr h}^\Psi (\overline\BBG[\gw_{n\,s}])$ is uniformly bounded 
the $g\circ \overline\BBG[\gw_{n\,s}]$ are uniformly integrable in $L^1_\gr(\Gw)$ by \rlemma{equiv}. Therefore the $(g\circ u_n)^+$
are uniformly integrable in $L^1_\gr(\Gw)$. Similarly 
\bel{Delta4'}g(x,u_n)\geq \gth\left(g(x,-v'_k)+g(x,-\overline\BBG(\gw^-_s)\right)-\ell
\ee
and the $(g\circ u_n)^-$ are also  uniformly integrable in $L^1_\gr(\Gw)$. The conclusion follows in the same way  as in \rth{delta2}.\qeda


\end {document}